\documentclass[a4paper]{article}

\usepackage{amsmath, amsfonts,amssymb,amsthm}

\usepackage{amssymb}
\newtheorem{theorem}{Theorem}[section]
\newtheorem{lemma}{Lemma}[section]
\newtheorem{proposition}{Proposition}[section]

\theoremstyle{definition}

\newtheorem{example}{Example}[section]

\theoremstyle{remark}
\newtheorem{remark}{Remark}[section]

\numberwithin{equation}{section}
\newcommand{\R}{\mathbb{R}}
\newcommand{\dist}{\mathrm{dist}}

\allowdisplaybreaks

\begin{document}
\title{Dynamic System of Neurons on a Complete Graph: Synchronization}

\author{S.A.~Pirogov\footnote{IITP RAS,HSE. The study has been funded by the Russian Science Foundation (project № 24-11-00123)} \and  A.N.~Rybko\footnote{MIPT} \and   D.D.~Pervushin\footnote{Skoltech} \and E.N.~Petrova\footnote{IITP RAS}}

\maketitle

\begin{abstract}
The net of $N$ ``physical’’ neurons is considered as a dynamical system. These neurons form a complete graph. The state of any neuron is its electric potential. The potential linearly increases until reaches its maximal value. Then it falls to zero and the neuron sends spikes to all other neurons. Having got a spike any neuron changes its state by some given function. We study the problem of a synchronization of the net. We glue the maximal value of potential to zero and so consider the state of any neuron as a point on the circle. We prove that under some conditions the states of neurons considered as points on the circle finally (when time turns to infinity) form a set which is time-invariant modulo its uniform rotation along the circle.
\end{abstract}

\section{Introduction}

We consider a symplified model of a neuron net. Each neuron is characterized by its potential $x$.  The value of $x$ is between $0$ and the maximal value which we assume to be equal to $1$. Neurons form the vertices of oriented graph $G(N,R)$ with $N$ vertices and $R$ edges. The potential $x$ of each neuron increses with time as $x(t)=x(0)+ct, \; c>0,$ until it reaches the value $1$. After that the value falls to $0$ and the neuron sends the spikes to its neubours on the graph. When a neuron at state $x$ gets a spike its state changes to $W(x)=x+f(x)\geq 0.$ The function $f(x)$ can have positive and negative values. We suppose that $W(x)$ can not exceed $1$. Also we suppose that $W(x)$ is monotone, i.e., if $x_1 \leq x_2$ then $W(x_1) \leq W(x_2).$ If several neurons reach the value $1$ simultaneously by continuous increasing of potentials or by jumps, then they act on their neighbours successively. Of course, the order of their actions is not important, and as a result the potentials of all these neurons take the value $0$. Then these neurons will have equal potentials forever (we use the term ``cluster'' for them).

Now we consider the case when the graph $G(N,R)$ is complete, i.e.\ each neuron is a neighbour of all the others. We consider the set of the potentials of $N$ neurons as a finite ordered sequence $0\leq x_1 \leq x_2 \leq \ldots \leq 1.$ When $x_N$ reaches the value $1$ the values of other neurons, $x_i$, are changed by jumps, but their order is conserved due to monotonicity, while $x_N$ falls to $0$. Thus the new value of $x_1$ is $0$ and the new values of the sequence are:
$$
\tilde x_1 =0, \; \tilde x_2=W(x_1), \; \tilde x_3=W(x_2), \ldots , \; \tilde x_N=W(x_{N-1})
$$
if $x_{N-1} < 1.$ If, conversily, $x_{N-1} = 1$, but $x_{N-2} < 1$ we have
$$
\tilde x_1 =0, \; \tilde x_2=0, \; \tilde x_3=W(x_1), \tilde x_4=W(x_2), \ldots , \, \tilde x_N=W(x_{N-2}).
$$
Considering all such cases we finish with the case when all $x_i$ equal $1$, thus $\tilde x_i=0$ for all $i$. When all jumps are accomplished, the values of $x_i$
continue  to increase linearly in time until the maximal of them reaches $1$. Then the corresponding neuron sends spikes to those neurons whose values have not reached $1$. It is natural to consider our system at discrete moments of time, namely, when one of the potentials reaches the value $1$. This discrete time dynamical system has a trivial fixed point
$x_1=x_2= \ldots =x_N=1$ when all $N$ neurons form one cluster. Later we shall look for non-trivial fixed points.

It is convenient to introduce a new variable $y=1-x$ and a new function $V(y)=1-W(x)\geq 0$. The monotonicity of $W$ is equivalent to the monotonicity of $V$. The fall of $x$ from $1$ to $0$ corresponds to the uplift of $y$ from $0$ to $1$; the linear increasing of $x$ with rate $c$ corresponds to the linear decreasing of $y$ with the same rate; the jump from $x$ to $W(x)$ corresponds to the jump from $y$ to $V(y)$. The uplift of $y$ takes place when $y$ reaches $0$. The trivial fixed point corresponds to $y_1=y_2=\ldots =y_N=0.$ Again we order the values of $y$ corresponding to $N$ neurons: $0 \leq y_1 \leq y_2 \leq y_N \leq 1.$

It is useful to think of the dynamic system $T(X)$ as a dynamic system on the circle of length one: there are $N$  particles  moving on the circle. 
Denote the coordinates of the particles by $\{ x_0(t), \ldots , x_{N-1}(t)\}$. There is a graph of connections between the particles, $G(N,R).$ 
And there is an influence function, $f(x)$, defined on $[0,1]$, such that  $-1 \leq f(x) \leq 1$ and $f(1)=0$.
We consider the origin of the circle as a special double point, one point with coordinate $0$, the other one with coordinate $1.$
The particles are moving clockwise with rate $1$ till one of them reaches the point $1.$ Immediately after that the coordinate of this particles changes to $0$, and the coordinate of any other particle changes according to the following rule. If there is no edge in the grapgh $G(N,R)$ connecting a particle with the  one 
that reached $1$, its coordinate is not changed. If there is an edge connecting a particle with the one that reached $1$, then its coordinate, $x(t)$, is changed to 
$x(t)+f(x(t))$. It means that this particle jumps, clockwise or countercleckwise, depending on the sign of $f(x(t))$. All the particles that a connected with the one that reached $1$ perform the jumps simultaneously. In this paper we suppose that $f(x)\geq 0$ and consider the mean field model, i.e., when the graph 
$G(N,R)$ is complete. We supose also that the particles cannot jump over each other, i.e., for $x>x',$ $f(x') - f(x) \leq x-x'.$ If several particles, say, 
with coordinates $x_1 \leq x_2(t) \leq \ldots \leq x_k$ jump to $1$ 
simultaneously (this happens if $f(x_1)=1-x_1, \ldots , f(x_k)=1-x_k$), then immediately all these particles induce the avalanche of jumps of their neighbours in the graph.

Let $t_1, \ldots , t_n \ldots$ be the moments when a particles reaches the point $1$.
We shall investigate the asymptotic behaviour of $T(X(t_n))=T(X(n))$ as $n\to \infty$. Since $f(1)=0,$ there is a trivial stationary state
$\bar X = \{ \bar x_0=1, \ldots , \bar x_i=1, \ldots , \bar x_{N-1}=1 \}$ when all  particles form a single cluster and this cluster rotates on the circle with rate $1$
without any jumps. Our goal is to to find all the stationary states of the dynamic system $T(X(n))$ and investigate the convergence to the stationary states.

We are considering a very rough mathematical model, omitting many important properties of the biological system in order to obtain rigorous mathematical results. We hope to obtain further results for a model that takes into account some of these properties. We think that such a class of dynamical systems is interesting also from mathematical point of view.

In recent research on synchronization phenomenon in physical and biological models the Kuramoto method is popular. By physical nmethod he found synchronization of two interacting oscillators. It seems that for real complicated physical and biological processes as, for example, those considered in \cite{1,2}, it is unlikely to obtain regorous mathematical results nowadays. The work \cite{3} contains some mathematical theorems for a model close to the one we are considering in this article. In \cite{4,5} the synchronization phenomenon has been proven for a large network of a different nature, but perhaps the ideas of this study may be useful for our future research.

\section{Discrete time dynamical system and fixed points}

Let the values of $y_i$ after the uplift and all the jumps be $0< y_1<y_2<\ldots <y_{N-1} < y_N=1.$
We call the sequence $Y=\{y_i\}$  {\em a configuration}, or {\em a particle}, or  {\em a point}.
The linear decreasing of $y_i$ and the uplift turns this configuration to $0< y'_1<y'_2<\ldots <y'_{N-1} < y'_N=1,$ where
\begin{align*}
  y'_1&=y_2-y_1,\\
        y'_2&=y_3-y_1,\\
\ldots &\ldots\ldots \\
              y'_{N-1}&=1-y_1.
\end{align*}
Further,  the jumps turn these values of $y'_i$ to
\begin{align*}
  \tilde y_1&=V(y_2-y_1),\\
        \tilde y_2&=V(y_3-y_1),\\
\ldots &\ldots\ldots\\
  \tilde y_{N-1}&=V(1-y_1),\\
                  \tilde y_N&=1.
\end{align*}

Now our aim is to look for a non-trivial fixed point of the map $Y \longmapsto \widetilde Y.$ Under some natural condition a fixed point exists and is unique. Firstly remark that if the function $V(y)$ is strictly increasing on the interval $(0,1)$ and $0< V(y) <1$ for $0<y<1$, then the map $Y \longmapsto \widetilde Y$
translates the open symplex $ \mathcal{Y}=\{ Y: 0< y_1<y_2<\ldots <y_{N-1} < 1\} $ into itself. Under additional condition of continuity we have the following result.
\begin{theorem}\label{th1}
  If the function $V(y)$ is continuous and strictly increasing on the interval $(0,1)$ and $0< V(y) <1$ for $0<y<1$, then the map $Y \longmapsto \widetilde Y$ has a unique fixed point in the open symplex $\mathcal{Y}$.
\end{theorem}

\begin{proof}
  Since $V$ is continuous and strictly increasing on $(0,1)$, $0< V(y) <1$ for $0<y<1$, by monotonicity we extend $V$ to the segment $[0,1]$, $0<\leq V(y) \leq 1$ for $0 \leq y \leq 1$ and $V$ is continuous on $[0,1]$. Consider the configuration 
  \begin{equation}\label{eq1}
    Y=  \{ 0< y_1<y_2<\ldots <y_{N-1} < 1\}
    \end{equation}
  in the open symplex $\mathcal{Y}$. The dynamicas of the system is defined as a composition of the affine map
  \begin{align}
    \label{eq2}
y'_{N-1}&=1-y_1, \nonumber \\
        y'_{N-2}&=y_{N-1}-y_1,\\
\ldots &\ldots\ldots \nonumber \\
              y'_{1}&=y_2 - y_1 \nonumber
  \end{align}
  and generally non-linear map that transforms a point $y$ to $V(y)$. Thus,
  \begin{align}
    \label{eq3}
\tilde y_{N-1}&=V(1-y_1),\nonumber \\
        \tilde y_{N-2}&=V(y_{N-1}-y_1),\\
 &\ldots \nonumber \\
              \tilde y_{1}&=V(y_2 - y_1) . \nonumber 
  \end{align}
  So we get a map $Y \longmapsto \tilde Y$ on the symplex $\mathcal{Y}$, which we can extend by continuity to the map on the closed symplex
  \begin{equation}
    \label{eq4}
\overline{\mathcal{Y}} = \{ Y: 0 \leq y_1 \leq y_2 \leq y_{N-1} \leq 1 \}.
\end{equation}
By the Bronwer fixed-point theorem the map $Y \longmapsto \tilde Y$ has a fixed point in $\overline{\mathcal{Y}} $.
We shall prove that this fixed point is contained in the open symplex $\mathcal{Y}$, i.e., the inequalities \eqref{eq4} are strict for it. The equations for a fixed point are:
\begin{align}\label{eq5}
y_{N-1}&=V(1-y_1),\nonumber \\
        y_{N-2}&=V(y_{N-1}-y_1),\\
 &\ldots \nonumber \\
              y_{1}&=V(y_2 - y_1). \nonumber 
\end{align}

a) If we suppose $y_1 =0$ then from \eqref{eq5} it follows that $V(y_2)=0$ and so $y_2=0$, and $V(0) =0$ by strict monotonicity. Again from \eqref{eq5} it follows
$V(y_3)=0$ and  $y_3=0$ and so on. Finally we have $y_{n-1}=0$, consequently, $0=V(1)$ which contradicts the strong monotonicity.

b) If we suppose $y_{N-1} =1$, then from \eqref{eq5} we see $1=V(1-y_1)$ and so $V(1)=1$ and $y_1=0$ by a strict monotonicity. But $y_1=0$ is impossible by a).

c) If we suppose $y_k=y_{k-1}$, then from \eqref{eq5} and a strict monotonicity of $V$ it follows that $y_{k+1}=y_{k}$, and so on. Finally we get $y_{N-1}=y_{N-2}$ and
$1=y_{N-1}$ which is impossible by b).

So the fixed point of the map $Y\longmapsto \tilde Y$ is contained in the open symplex $\mathcal{Y}$.
Let us denote this point by $Y^*$ and prove that it is unique.  Suppose that $\hat Y$ is a fixed point of the map $Y\longmapsto \tilde Y$ and let us prove that
$\hat Y = Y^*$. It is sufficient to prove that
$\hat y_1 = y_1^*$ because any solution of \eqref{eq5} is completely determined by $y_1$ due to the strict monotonicity of $V$. Suppose that
$\hat y_1 > y_1^*$. Then from \eqref{eq5} it follows $\hat y_{N-1} < y_{N-1}^*$. Then  $\hat y_{N-2} < y_{N-2}^*$ due to a strict monotonicity of $V$, and so on, until
$\hat y_{2} < y_{2}^*$ and $\hat y_{1} < y_{1}^*$ which contradicts our hypothesis. The possibility of $\hat y_{1} < y_{1}^*$ leads to a contradiction by the same line of arguments. We conclude that the fixed point of the map $Y\longmapsto \tilde Y$ is unique.
\end{proof}

\begin{example}
  Let $N=2$. Then the configuration $Y$ consists of one point, $0<y_1<1$ and the map  $Y \longmapsto \tilde Y$ is simply
  \begin{equation} \label{eq6}
    \tilde y_1 = V(1-y_1).
  \end{equation}
  The function $y_1-V(1-y_1)$ is strictly monotone and takes values $-V(1) <0$ and $1-V(0) >0$ at the ends of the segment $[0,1]$. Thus the equation
  \begin{equation} \label{eq7}
    y=V(1-y)
    \end{equation}
    has a unique solution $y^*$.

    A new phenomenon appears when $V(0)=0$ and $V(1)=1$. Consider the iteration of the map $Y \longmapsto \tilde Y$:
    \begin{equation}\label{eq8}
      \widetilde{\tilde y_1} = V(1-V(1-y_1)).
      \end{equation}
      The map \eqref{eq8} has two additional solutions apart of $y^*_1$, namely, $y=0$ and $y=1$. If $V$ is differentiable, $V'(0)=c_0,$ $V'(1)=c_1$, then for small $y_1$ we have
      \begin{equation}\label{eq9}
        \widetilde{\tilde y_1} = c_0 \bigl(1-(1-c_1y_1)\bigr) + o(y_1) = c_0c_1y_1 + o(y_1).
      \end{equation}
      In the case when $c_0 c_1 <1$ the even iterations of $Y \longmapsto \tilde Y$ applied to $y_1$ tend to $0$ and the odd iterations tend to $1$. It means that the map $Y \longmapsto \tilde Y$ has an attracting periodic points $\tilde 0=1,$ $\tilde 1 =0$.
    \end{example}

    \begin{example}
      Let $V(y)=ay, \; a\leq 1.$ Then we have:
      \begin{align}\label{eq11}
   \tilde y_{n-1}&=a(1-y_1),\nonumber \\
        \tilde y_{n-2}&=a(y_{N-1}-y_1),\\
 &\ldots \nonumber \\
              \tilde y_{1}&=a(y_2 - y_1) . \nonumber 
  \end{align}     
  Let us introduce new variables:
  \begin{align}\label{eq12}
    z_1&=y_1, \nonumber\\
    z_2&=y_2-y_1, \nonumber\\
    z_3&=y_3-y_2\\
       &\ldots  \nonumber\\
    z_N&=1-y_{N-1}.  \nonumber
  \end{align}
  Then we have $\sum_1^N z_i=1$ and
   \begin{align}\label{eq13}
   &\tilde z_{1} =az_2,\nonumber \\
        &\tilde z_{2}=az_3,\\
 &\ldots \nonumber \\
     &\tilde z_{N-1}=az_N,\nonumber \\
     &\tilde z_N = (1-a) +az_1              . \nonumber 
   \end{align}
   The affine map defined by \eqref{eq13} has a fixed point $z^*_N=A, \; z^*_{N-1} =aA, \ldots ,$ $z^*_1=a^{N-1}A,$ where $A=\bigl(\sum_1^{N-1} a^k \bigr)^{-1}$. This fixed point is attractive for $a<1$ and neutral for $a=1$.
 \end{example}

 \begin{example}
 Physically, the mapping $W$ (in variables $x$) or $V$ (in variables $y=1-x$) is applied when one of the neurons reaches its maximum potential, $x=1$, and
 then the potential falls to $x=0$. At this moment the neuron sends the spikes to other neurons. We can consider the  situation where different neurons have different functions $W$. So for $N$ neurons we have $N$ different functions $W$ (or different functions $V$). Let us denote them by $V_1, \ldots , V_N$. Let us consider a linear case
 \begin{equation}\label{eq14}
   V_k(y)=a_ky, \; a_k \leq 1.
 \end{equation}
 Denote by $T$ the operator in $\R^N$ acting as follows:
 \begin{equation}\label{eq15}
   T(z_1, \ldots, z_n) =(z_2, \ldots , z_N, z_1).
   \end{equation}
   Denote by $e$ the vector $(0,0,\ldots , 0,1)$, and we denote $(1-a)$ by $\bar a$. Then the affine map \eqref{eq13} can be written as
   \begin{equation}\label{eq16}
     B_az=aTz +\bar a e.
   \end{equation}
   Now we express an operator corresponding to a consecutive discharge of all $N$ neurons:
   \begin{equation}\label{eq17}
     B=B_{a_N} B_{a_{N-1}} \ldots B_{a_1}.
     \end{equation}
     Starting calculations from
     \begin{align}\label{eq18}
       B_{a_2} B_{a_{1}}z&=a_2(a_1T^2z+\bar a_1Te)+a_2e=\nonumber \\
                         &=a_2a_1T^2z+a_2 \bar a_1Te+\bar a_2 e
     \end{align}
     we finally get
   \begin{align}\label{eq19}
     B_{a_N} B_{a_{N-1}}\ldots B_{a_{1}}z&=a_Na_{N-1} \ldots a_1 T^N z+\nonumber \\
                                         &\quad {} + a_N\ldots a_2\bar a_1T^{N-1}e + \ldots +a_N \bar a_{N-1}Te                                           
                         +\bar a_N e.
     \end{align}  
     Since $T^Nz=z,$ the operator $B$ is simply
     \begin{equation}\label{eq20}
     Bz-a_N\ldots a_1 z+g,
     \end{equation}
     where
     \begin{equation}\label{eq21}
       g=(a_N \ldots a_2 \bar a_1, a_N \ldots \bar a_2, \ldots, a_n\bar a_{N-1}, \bar a_N).
     \end{equation}
     Thus  $B$ is contracting operator (if not all of $a_k$ are equal to $1$) and so it has a unique attracting fixed point
     \begin{equation}\label{eq22}
       z^*=g/(1-a_N\ldots a_1),
     \end{equation}
     and $z^*$ is a global attractor of our dynamical system.
   \end{example}
   
     We now return to the case when all neurons are identical and
     $$
     V(y)=ay, \; a \leq 1.
     $$
     Let us consider the empirical distribution $\mu^{(N)}$ corresponding to the fixed point $Y^*$, i.e., the distribution assigning the measure $1/N$ to each $y_i$.
     Let $a$ depend on $N$. Three different types of limit behavior of  $\mu^{(N)}$ are possible:
\begin{itemize}
    \item[1)] $1-a=o(1/N), \mbox{ i.e., } (1-a)N\to 0.$ Then  $\mu^{(N)}$ tends to uniform destribution on $[0,1]$.
    \item[2)]  $(1-a)N\to \alpha >0.$
      Then  $\mu^{(N)}$ tends to a distribution on $[0,1]$ having the density $1/(\alpha y + \beta),$ where $\beta$ is defined from the condition
      \begin{equation}\label{eq23}
        \int\limits_0^1 \frac{{\rm d} y}{\alpha y + \beta} =1.
        \end{equation}
      \item[3)]
        $(1-a)N \to \infty.$ Then $\mu^{(N)}$ tends to an atom at zero point, i.e.\ to a $\delta$-measure.
      \end{itemize}
      These results are obtained by easy direct calculations.

      \begin{remark}
        The map $Y\longmapsto \tilde Y$ is a composition of two maps: affine shift and a non-linear map $y_i \longmapsto V(y_i)$ for all $y_i$. Of course we can consider an opposite composition which turns the configuration
        $$
0< y_1 <y_2, \ldots y_{N-1}<1
        $$
        to
        $$
0< V(y_1) <V(y_2), \ldots V(y_{N-1}) <1
$$
and then to
\begin{align*}
  \hat y_{1}&=V(y_2) - V(y_1) ,\\
  &\ldots \\
 \hat y_{N-2}&=V(y_{N-1}), \\
        \hat y_{N-1}&=1-V(y_{1}).
   \end{align*}
   Obviously, the fixed point of the map    $Y\longmapsto \hat Y$ can be obtained from the fixed point of the map    $Y\longmapsto \tilde Y$
   by applying pointwise the function $V^{-1}$, which is the inverse to $V$.
 \end{remark}
 
\section{Trapezoidal $\boldsymbol{f(x)}$ dynamics}

Now let us consider the case when the functions $w$ are continuous but not strictly monotone. Remind that $W(x)=x+f(x).$ Now the function $f(x)$ depends on parameter $h$, $0<h<1.$ Namely,
$$
f_h(x)=\begin{cases}
  h & \mbox{for } x \in (0,1-h],\\
  1-x & \mbox{for } x\in [1-h, 1).
  \end{cases}
  $$
  It means that
  $$
W(x)=\begin{cases}
  x+h & \mbox{for } x \in (0,1-h],\\
  1 & \mbox{for } x\in [1-h, 1).
  \end{cases}
  $$  
  For any initial configuration $0<x_1 < x_2 < \ldots < x_{N-1} <1$, after applying $W$ all the points $x_i$ in the interval $[1-h,1)$ jump to $1$. 
  The other particles are shifted by $h$ in positive direction. The new particles that appeared in $1$ do the same. After repeating this process, as a result we get a cluster of particles concentrated in $1$, which then fall to $0$, and there are no particles in $[1-h, 1)$.
  Further, after linear increasing of coordinates, at the moment when any particle appears in $1$ the process starts again. For this system the following theorem takes place.
  \begin{theorem}
    \label{th2}
    For the function $W(x)=x+ f_h (x)$ the following statements hold true.
    \begin{itemize}
    \item[\rm {1)}] If $Nh \geq 2,$ then for any initial configuration $X$ the trajectory of the dynamical system under consideration reaches the trivial fixed point: all particles will be glued in one cluster.
    \item[\rm {2)}]
If $Nh<1$, then for any divisor $k$ of $N$ there are fixed points which consist of the regularly placed clusters of $N/k$ particles, and any trajectory reaches one of these fixed points.
\item[\rm {3)}]
  If $1<Nh <2$ and $k$ is a divisor of $N$ such that $Nh - \frac{Nh}{k} < 1$, then there are fixed points
$\bar X=\{ Z_1, \ldots , Z_k \}$
  which consist of regularly placed $N/k$ clusters of $k$ particles, and any trajectory reaches one of these fixed points. If for any divisor $k$, $Nh-\frac{Nh}{k} \geq 1$, then, like in the case {\rm 1)}, any initial state converges to a trivial stationary state, i.e., to one cluster.
  \end{itemize}
    \end{theorem}

    To prove Theorem~\ref{th2} we need the following lemma.
    \begin{lemma}
      \label{l1}
      Let the state of the dynamical system $T(X(t))$ consist of $k$ clusters, $X(t)=\{Z_1(t), \ldots , Z_k(t)\}$.
Denote by $n_i$ the number of particles in the cluster $Z_i.$ Thus, a  cluster $Z_i$ can be described by a pair $\{n_i,x_i\}$, the number of particles and their common coordinate. If there are
      two clusters, say, $Z_i(t)$ and $ Z_l \}$ with different number of particles,
      $n_i \neq n_l$, then after a finite time interval one of the  clusters will increase its length.
      \end{lemma}
      \begin{proof}
       Consider the clusters $Z_0$ and $Z_1$. The distance between the clusters is not changed between the moments of jumps. At the moment of the jump (corresponding to $x_1=1$) the distance between the clusters, $|x_1 - x_0$, changes to $|x'_1 -x'_{0}|$ as follows:
      \begin{align}
        \label{eq24}
        &|x'_1 -x'_0|=\begin{cases}
          |x_1 -x_0|  & \mbox{if } |x_1 -x_0| > n_1h,\\
          0 & \mbox{if } |x_1 -x_0| \leq n_1h.
          \end{cases}\\
        \label{eq25}
        &|x'_2 -x'_1|=\begin{cases}
          |x_1 -x_0|  & \mbox{if } |x'_2 -x'_1| > 0,\\
          0 & \mbox{if } |x'_1 - x'_0| = 0 \mbox{ and } (n_0 +  n_1)h  <  |x_2 -  x_0|,
          \end{cases}
        \end{align}
        and so on.  We see that for one rotation a small cluster jumps at a distance larger than a large cluster does.
        Thus if we assume that for an arbitrary number of rotations there is no sewing together for any clusters, we will find that all the clusters must have the same mass: $n_i \equiv n_j$, since for one rotation we have
  \begin{equation}
    \label{eq26}
    |x'_i -x'_j|= |x_i -x_j|  -h|n_i - n_j|,
  \end{equation}
  i.e., a large cluster pushes a small one stronger than the other way around.
\end{proof}

We now proceed to the proof of Theorem~\ref{th2}.
\begin{proof}

%\subsection*{Proof of Theorem~\ref{th2}}

We are going to describe all the configurations containing  clusters of equal mass which have not coupled. Let us consider the case  when $Nh<2$, one cluster
is in the state $x_0(t)=1$ and all other clusters contain $N/k$ particles. Let us enumerate all the distances between neighbouring clusters counter-clockwise.
The distances $ \dist (Z_m,Z_{m+1})$ between neighbouring clusters must be greater than $hN/k$ for $m>0$, since they do not change until  the cluster
$Z_m$ reaches the point $1$, but if at the moment before the jump the distance becomes less than or equal to $hN/k$,  this pair of clusters couples. Immediately after any jump the distance $\dist (Z_0,Z_1)$ must be greater than $0$. It is easy to verify that for our dynamics
 any configuration that consists of $k$ clusters, $X=\{ Z_0, \ldots , Z_{k-1}\}$ such that $\dist(Z_0,Z_1)=|x_1 - x_0|>0$, $\dist(Z_{m+1}, Z_m)=|x_{m+1} - x_m|>Nh/k$
 at any moment immediately after the jumps is stable up to permutation of coordinates:
$\{x_0=1, \ldots , x_{k-1}\} \to \{x_1=1, \ldots , x_{k-1}, x_0\},$
 %any vector that satiafies the conditions:
%$$
%\dist (X_0, X_1) >0, \quad \dist(X_m,X_{m+1})> \frac{Nh}{k}
%$$
%
%$$\{x_0=0, \ldots , x_{k-1}\} \to \{x_1=0, \ldots , x_{k-2}, x_0\},$$
and this proves items 2) and 3) of the theorem.

If
\begin{equation}\label{eq27}
  \sum_{i=0}^{k-1} \dist (Z_i, Z_{i+1})=\Bigl( N-\frac{N}{k}\Bigr) h \geq 1,
\end{equation}
then it is impossible to put the clusters on the circle of length $1$ with all the distances $\dist (Z_m,Z_{m+1})> \frac{N}{k} h$. This proves point 1) of the theorem.
\end{proof}

It is easy to see that if $hN\geq 2$, the inequalities \eqref{eq27} are never satisfied, except the case when $k=1$ which corresponds to a trivial stationary state when all the particles have coupled in one cluster. All stationary states with equal  divisors $k$ of $N$ have positive Lebesgue measure -- the position of clusters on the circle may be moved under inequalities \eqref{eq27}. Nevertheless it is natural to conject that if $N$ is large enough and $Nh$ is of order of a constant, then the trivial stationary state is dominant, i.e.,  for ``typical'' initial configurations there is onvergence to this trivial stationary state. To confirm this conjecture, we will show a few simple examples.

Let $N$ be a large prime number, then we have only two possibilities: either $k=1$ or $k=N$. Let us take as initial state the uniform distribution of $N$ independent particles on the circle. This random configuration converges as $N\to \infty$ to a configuration of Poisson random field with intens
ity $N$ on the unit circle.

%%%%%%%%%%%%%%
If at least for one pair of neighbouring particles, $x_m$ and $x_{m+1}$, the distance $\dist (x_m, x_{m+1})$ is less than $h$, then they will couple in a common cluster when the particle $m$ comes to point $1$ on the circle, and according to Theorem~\ref{th2} all particles after a finite time will couple to a common cluster
(trivial stationary state). But for Poisson field the probability for all the distances between neighbouring particles to be greater than $h$ exponentially decreases with $N$.

\begin{proposition}\label{p2}
  Consider a model consisting of different particles, where each particle $i$ sends the impulse $h+\xi_i$ to other particles, $i=1, \ldots , N$. Assume that each $\xi_i$ is a small random variable with a positive continuous density defined on an arbitrarily small interval around zero and all the $\xi_i$ are mutually independent.
  Then for any initial state, with probability $1$ within a finite interval of time (that could be large if $i$ is small), the trajectory of the dynamical system converges to the trivial stationary state. 
\end{proposition}
\begin{proof}
  Let us define the ``force'' of jumps of particles when a cluster $X_m =\{i_1, \ldots , i_m \}$ comes to point $1$. The jump of any particle $k$ from $X_m$ is equal to
  \begin{equation}\label{i}
    F_k=\sum_{i_k\in X_m} h+\xi(i_k).
  \end{equation}
  It is clear that for any pair of clusters, $X_m$ and $X_l$, the probability that $F_k=F_l$ is equal to zero. Applying Lemma~\ref{l1}, by repeating its arguments we immediately have that for any initial configuration the system converges to trivial stationary state with probability one.
\end{proof}

\begin{proposition}\label{p3}
 % Let the conditions 41-44 be fulfilled.
  Let $\bar f(x)$ be continuous, $\bar f(x) \geq f_h (x)$ and $\bar f(x) = f_h (x)$ for $x \in [1-h,1]$, and let for  any $x_1> x_2$, $\bar f(x_1) \leq f(x-2)$.
  If $Nh \leq 2$, then for any initial state the dynamical system with function $\bar f(x)$ converges to trivial stationary state.
\end{proposition}
\begin{proof}
  We prove by contradiction: it is clear that at any moment there exists minimal cluster $Z_i(t)$.
  Its mass, $n(Z_i(t)),$ cannot be larger than $N/2$, hence,
   $hn(Z_i(t)) \leq 1.$
   Then  $\sum_{j\neq i} h n(Z_j) \geq 1$. Let us wait till the cluster $Z_i(t)$ be in the last position on the circle.
   It is easy to see that among the  other clusters at this moment of time there exists a cluster $Z_m$ such that $\dist (Z_m, Z_{m+1}) \leq h n(Z_m)$ and
   % at the moment $t$
   when
  the cluster $Z_m$ comes to the point $1$ the cluster $Z_{m+1}$ couples with $Z_m$ (as $\dist (Z_m, Z_{m+1})$ have not increased upon time, because
  $\bar f(x)$ is nonincreasing).
\end{proof}

\section{General case of monotone function $\boldsymbol{f(x)}$ on the circle}

Let us consider the generalization of the previous section, namely:
\begin{itemize}
\item[{\rm (i)}]
  let $f(x)$ be an arbitrary continuous decreasing function on the circle which is between two linear functions:
%\begin{align}\label{eq28}
  $1-x > f(x) > \varepsilon - \varepsilon x,$;
  \item[{\rm (ii)}]
    $-1+\varepsilon < f'(x) < -\varepsilon$;
  \item[{\rm (iii)}]
    $f''$ is bounded on $[0,1]$;
  \item[{\rm (iv)}]   $f(1)=0$.
\end{itemize}
    The next theorem is a generalization of Theorem~\ref{th2}. 
\begin{theorem}\label{th3}
  Let $f(x)$ satisfy conditions {\rm (i)} -- {\rm (iv)}. Then
  \begin{itemize}
  \item[1)]
    The dynamic system $T_N(x)$ of $N$ particles at discrete moments of time, $\delta t,$ -- at the moments when recurrent particle comes to point $1$ -- has one stationary state, $\bar X=\{ \bar x_0=1, \bar x_z < \bar x_0, \ldots , \bar x_{i+1}<\bar x_i, \ldots \bar x_{N-1}\}$ for arbitrary initial configuration of isolated particles.
  \item[2)]
    For any initial state of $N$ isolated particles the dynamic system converges to this stationary state with exponential rate in $l_\infty $ norm.
    \end{itemize}
  \end{theorem}
  \begin{proof}
   According to Theorem~\ref{th1} the dynamic system has one stationary state $\bar x = \{ \bar x_1, \ldots , \bar x_{N}$, $\bar x_1=1,$ $\bar x_{n-1} < \bar x_n, \ldots  < \bar x_N <0$. We shall prove now that for any initial state $\{ x_0, x_1, \ldots , x_N\}$ of isolated particles we have convergence to $\bar x_N$ with exponential rate:
                \begin{equation}\label{eq35}
                  \lim_{n\to\infty} T^n (x) = \bar x.
                  \end{equation}
                  The evolution defined by $T(x)$ is a product of two operators: the operator $A$ defines the jumps of the particles, and the operator $B$ defines their shift and renumeration.
                  %%%%%%%%%%%%%%%
  The operators $A$ and $B$ do not commute. We take $T'=AB$ and $T''=BA$, then for a stable point $\bar X'$ for $T'$ and a stable point
  $X''$ for $T''$ we have trivial relations: $\bar X''=B(\bar X')$ and $\bar X'=A(X'')$. Thus we can use either $T'$ or $T''$ as $T$ when it is convenient for the proof of the theorem.

 Let us write down the equations for $T(x)$ and for a fixed point $\bar x$.
    \begin{align}\label{eq29}
      T(x)&=T(x_0=1;x_1;\ldots ;x_{N-1})=BA(x),\nonumber \\
      B(x)&=\{y'_0=1, \ldots, y'_{N-1}\},
      \intertext{where}
      y'_2&=x_2+(1-x_1), \, \ldots \, , \, y'_i=x_{i+1}+(1-x_1), \, \ldots \,  , \, y'_{N-1}=1-x_1. \nonumber 
    \end{align}
    \begin{align}\label{eq30}
      A(B(x))&=\{y''_0=1, y''_1, \ldots , y''_i, \ldots , y_{N-1} \}, 
      \intertext{where}
      y''_1&=y'_1+ f(y'_1), \, \ldots , \quad y''_i=y'_i+f(y'_i), \, \ldots , \nonumber \\
      y''_{N-2}&=y'_{N-2} + f(y'_{N-2}), \;  y''_{N-1}=(1-x_1) +f(1-x_1) . \nonumber 
        \end{align}
        For a fixed point $\bar x$ we have:
        \begin{align}\label{eq31}
          &AB(\bar x)=\bar x = \{\bar x_0, \ldots , \bar x_i, \ldots , \bar x_{N-1} \}\\
            \intertext{where}
          \bar x_{N-1} &= (1-\bar x_1) + f(1-\bar x_1), \nonumber \\
          \bar x_{i} &= \bar x_{i+1} + (1-\bar x_i) + f(\bar x_{i+1} +(1-\bar x_i)), \ldots , \nonumber \\
            \bar x_{1} &= \bar x_{2} + (1-\bar x_1) + f(\bar x_{2} +(1-\bar x_1))  .     \label{eq32}
            \end{align} 
    It is easy to see that $\bar x$ is the unique solution of \eqref{eq31}--\eqref{eq32}, $\bar x = \{x_1=0, \ldots , \bar x_i > \bar x_{i+1}, \ldots , \bar x_{N-1}>0 \}$
                because of monotonicity of the function $f(x)$.
                Note that because of monotonicity of $f(x)$, after two rotations all configurations $T^{2N}(x)$ acquire the following property: uniformly on $x$,
                \begin{align}
                  &T^{2N}(x_{N-1}) >c_1>0, \label{eq33}\\
                  &T^{2N}(x_1) <1-c_2>0 \label{eq34}
                \end{align}
                Thus, we can forget about the discontinuity of $f(x)$ in $x=1$ and consider $T(x)$ in the interval $[c_2, 1]$ with continuous function $f(x)$.
                In biological applications typically the function $g(x)$ is continuous on $[0,1]$, $g(x)=0$ at $1$, it increases on a short interval $[0,\varepsilon]$ and decreases on $[\varepsilon , 1]$, so the behaviour of our $f(x)$ and $g(x)$ will be the same for $T^k(x)$ for $k>2N$, if $c_1,>c$.            

                Mow we prove a simple useful lemma
                \begin{lemma}\label{l2}
                  We shall take $w(x)=BA(x)$. There exists $\delta >0$ such that for any state $x$ of isolated particles there exists its neighbourhood $u(x)$ such that for any $y \in u(x)$
                  \begin{equation}\label{eq350}
                    \| T(x) - T(y) \|_{l_\infty} < \|x-y\|_{l_\infty} (1-\delta).
                  \end{equation}
\end{lemma}
\begin{proof}
  By definition, $\|x-y\|_{l_\infty}$ is the maximum of the arches $l_i=\dist(x_i, y_i)$ on the circle,
\begin{equation}\label{eq36}
  \|x-y\|_{l_\infty} = \max \{ |x_i-y_i|, i=1, \ldots , N-1\}, \quad x_0=y_0=1.
  \end{equation}
  To estimate $\| T(x) - T(y) \|$, remind that $w=BA$ where $A$ defines the jumps and $B$ defines the shift. For small $l_i$, as $f''(x)$ is bounded, we have the estimation
  \begin{equation}\label{eq37}
    \sum_{i=1}^{N-1} \| A(x) - A(y) \| = \max \bigl | (x_i - y_i) f'(x_i) + o (x_i-y_i)\bigr |,
  \end{equation}
  where $-1 + \varepsilon < f'(x) < - \varepsilon$. %for sufficiently small neighbourhood of $x$.
  So
  \begin{equation}
    \label{eq38}
    \| A(x) - A(y) \| < \max_{1\leq i \leq N-1} \bigl | (x_i - y_i)  (1-\varepsilon) \bigr |
  \end{equation}
for all $Y$ in sufficiently small neighbourhood of $X$.
We may rewright \eqref{eq38} as follows:
 \begin{equation}
    \label{eq39}
    \| A(X) - A(Y) \| < \max\Bigl[ \max_{2\leq i \leq N-1} \Bigl [ \bigl | x_i - y_i \bigr |  \Bigl(1-\frac{\varepsilon}{2}\Bigr)\Bigr ] ,  \Bigl(1-\frac{\varepsilon}{2}\Bigr)
    \bigl |x_1 - y_1 \bigr | \Bigr ]
  \end{equation}

  Now let us estimate $\| BA(X) - BA(Y) \| .$ The map $ BA(X)$ shifts all coordinates $i=1, \ldots , N-2$ of vector  $ A(X)$ on $A(x_1)$ and renumbers them
  (note that renumbering does not change the norm of $ BA(X).$ Respectively, the map $B(A(Y))$ shifts the coordinates of vector $A(Y)$, and
 \begin{align}
    \label{eq40}
   &\| BA(X) - BA(Y) \| \\
   &\quad \leq \max\Bigl[ \max_{2\leq i \leq N-2} \Bigl [ \bigl | x_i - y_i \bigr |  \Bigl(1-\frac{\varepsilon}{2}\Bigr)
    + \Bigl(1-\frac{\varepsilon}{2}\Bigr)
    \bigl |x_1 - y_1 \bigr | \Bigr ] ,
    \Bigl(1-\frac{\varepsilon}{2}\Bigr)
    \bigl |x_1 - y_1 \bigr | \Bigr ].\nonumber 
  \end{align}
  It is easy to find now a small neighbourhood of $X$ which is contracting under $T(X)$, namely, if $\max_{2\leq i \leq N-2} (x_i - y_i )=0,$ then it is sufficient to take
  $$
  A \bigl |x_1 - y_1 \bigr | <  \bigl |x_1 - y_1 \bigr | \Bigl(1-\frac{\varepsilon}{2}\Bigr).
  $$
  If $\max_{2\leq i \leq N-2} (x_i - y_i ) \Bigl(1-\frac{\varepsilon}{2}\Bigr) > 0,$ then it is sufficient to take $ \bigl |x_1 - y_1 \bigr | \Bigl(1-\frac{\varepsilon}{2}\Bigr)$
  of order
   \begin{equation}
    \label{eq41}
\varepsilon^2 \max_{2\leq i \leq N-2} (x_i - y_i ) \Bigl(1-\frac{\varepsilon}{2}\Bigr) .
\end{equation}
We have that for the small neighbourhood of $X$, described above (let us denote it by $d(X)$),
\begin{equation}
  \label{eq42}
  \| BA(X) - BA(Y) \| \leq (1- \varepsilon) \| X- Y \|.
\end{equation}
Thus the lemma is proved.
\end{proof}

Now let us consider a segment $[X,\bar X],$ where $\bar X$ is such that $\{ \bar x_0=1, ber x_1 < \bar x_0, \ldots , \bar x_i < \bar x_{i-1}, \ldots , \bar x_{N-1} \geq 0.$ It is clear that all the points of this segment are isolated particles with $x_0=1.$ 
    
Consider a covering $S$ of the segment $[X,\bar X]$ by open intervals such that for any $\tilde X \in [X,\bar X]$ there is an interval $d(X) \in S$ such thst for any $Y \in d(X)$,
\begin{equation}
  \label{eq43}
  \| T(X) - T(Y) \| \leq (1- \varepsilon') \| X- Y \|.
\end{equation}
Since $[X,\bar X]$ is compact, we can choose a finite covering $\{d_1, \ldots , d_k$ by the intervals which satisfy \eqref{eq43} and such that
$\sum_{j=1}^k |d_j| -\| X- \bar X \|$ is arbitrary small, where $|d_j| $ is the length of the interval $d_j.$ 
Let us take $k$ points, $X_1, \ldots , X_k \in \| X- \bar X \|$ and the corresponding intervals $d(X_1), \ldots , d(X_k)$, and for any $j$ let us choose a point
$Z_j \in d(X_i) \cap d(X_{j+1}$. Consider a broken line $L$ of connected segments
\begin{align}
  &\big[T(X_1=X, T(Z_1 \in d(X_1) \cap d(X_2))\big] \cup \big[T(Z_1),T(X_2)\big] \cup \ldots \nonumber \\
  &\; \cup \big[T(X_j), T(Z_j \in d(X_j) \cap d(X_{j+1})) \big] \ldots\nonumber \\
  &\; \cup \big[T\big(Z_{k-1} \in d(X_{k-1}) \cap d(\bar X)\big), T(\bar X)\big]. \label{eq44}
\end{align}
By construction, the total length $d(L)$ of this broken line is not greater than $(1- \varepsilon') \| X- \bar X \|.$ Thus, 
$\| T(X) - \bar X \| \leq (1- \varepsilon') \| X- \bar X \|.$
Repeating further this procedure, namely, connecting by the segment $[ T(X) - \bar X ] $ and using the corresponding covering, we see that $X$ converges to $\bar X$ exponentially fast. The theorem is proved.
\end{proof}

Let us now formulate a simple proposition which is useful in similar situations.
\begin{proposition}
  \label{p4}
  Let a phase space $\mathcal{ X}$ of the discrete time dynamical system $T(X)$ in the metric space $\mathcal{Y}$ contain an invariant set $\tilde{\cal X}$ 
  which has exactly one stable point $\bar X.$ Suppose that for any $X \in \tilde{\cal X }$ there exists a neighbourhood $d(X)$ such that for any $Y \in d(X)$,
 \begin{equation}
   \label{eq45}
   \dist (X,Y) > \dist (T(X), T(Y)) (1-\varepsilon).
 \end{equation}
 Suppose that for any $X \in {\cal X}$, $[X, \bar X] \in \tilde {\cal X}.$ By $[X, \bar X] $ we mean a connected compact, its length is defined as the minimal length between $X$ and $\bar X$ in the metrical space ${\cal Y}.$  Then for any $X \in \tilde{\cal X}$, 
\begin{equation}
  \label{eq46}
  \lim_{n\to\infty} T^n(X) =\bar X
\end{equation}
exponentially fast. Also, if $\mathcal{Y}$ is a compact itself, then, due to \eqref{eq45}, $T(Y)$ has a unique fixed point $\bar X .$
\end{proposition}

The proof of this proposition is the same as the proof of Theorem \ref{th3}.

\end{document}